\documentclass[prl, sd, amsmath, amssymb, reprint]{revtex4-1}

\usepackage{amsmath}
\usepackage{amsfonts}
\usepackage{amssymb}
\usepackage{amsthm}
\usepackage{fancyhdr}
\usepackage{enumerate}
\usepackage[top=1in, bottom=1in, left=1in, right=1in]{geometry}
\usepackage{graphicx}
\usepackage[labelfont=bf]{caption}

\usepackage{caption}
\usepackage{subcaption}

\usepackage{verbatim}

\usepackage{wasysym}
\usepackage{hyperref}

\usepackage{wrapfig}

\begin{document}

\title{Comparing the Locking Threshold for Rings and Chains of Oscillators}
\author{Bertrand Ottino-L\"{o}ffler and Steven H. Strogatz}
\affiliation{Center for Applied Mathematics, Cornell University, Ithaca, New York 14853}
\date{\today}

\begin{abstract}
We present a case study of how topology can affect synchronization.  Specifically, we consider arrays of phase oscillators coupled in a ring or a chain topology.  Each ring is perfectly matched to a chain with the same initial conditions and the same random natural frequencies. The only difference is their boundary conditions: periodic for a ring, and open for a chain. For both topologies, stable phase-locked states exist if and only if the spread or ``width" of the natural frequencies is smaller than a critical value called the locking threshold (which depends on the boundary conditions and the particular realization of the frequencies). The central question is whether a ring synchronizes more readily than a chain.  We show that it usually does, but not always.  Rigorous bounds are derived for the ratio between the locking thresholds of a ring and its matched chain, for a variant of the Kuramoto model that also includes a wider family of models. 
\end{abstract}

\maketitle

\section{Introduction}

The Kuramoto model has been used to study the dynamics of synchronization in a wide variety of physical, chemical, and biological systems~\cite{kuramoto75, kuramoto84, strogatz00, pikovsky03, strogatz03, acebron05, dorfler14, pikovsky15, rodrigues16}. The model's governing equations can be written (unconventionally, but most usefully for our purposes) in the following dimensionless form: 
\begin{alignat}{1} \label{GeneralKura}
\dot\theta_k = \Gamma \eta_k + \sum_{j\in\mathcal{N}(k)}   \sin(\theta_j - \theta_k),
\end{alignat} 
for $k = 1, \ldots, N$, where $\theta_k$ is the phase of oscillator $k$, and the sum is over all of $k$'s neighbors $\mathcal{N}(k)$, as determined by the coupling graph.  By rescaling time in Eq.~\eqref{GeneralKura}, we have normalized the coupling strength to unity without loss of generality. The term $\Gamma \eta_k$ can then be interpreted as the scaled natural frequency of oscillator $k$. 

The motivation for this unusual notation is that we are going to regard $\eta = (\eta_1, \ldots, \eta_N)$ as a fixed frequency vector and $\Gamma \ge 0$ as an adjustable parameter controlling the spread of the natural frequencies. For instance, the components of $\eta$ could be chosen independently at random from a prescribed probability distribution. Then increasing $\Gamma$ would allow us to increase the ``width" of the set of frequencies $\left\{\Gamma \eta_1, \ldots,\Gamma \eta_N \right\}$. We will occasionally write $\omega_k := \Gamma \eta_k$ for brevity. 

In the simple case where $\Gamma=0$ and all the oscillators have $\omega_k=0$, the model has a stable fixed state with $\theta_k = 0$ for all $k$, for a broad class of coupling graphs. Now imagine increasing $\Gamma$ slightly to produce some variation among the $\omega_k$. Starting from an initial condition $\theta_k(0)=0$ and assuming a sufficiently small but nonzero $\Gamma$, the system will asymptotically approach a stable periodic solution of Eq.~\eqref{GeneralKura} in which all the oscillators run at the same constant frequency $\dot\theta_k \equiv \Omega$ for all $k$, for some constant $\Omega$. We call such a solution a stable phase-locked state.  But when $\Gamma$ gets too large, the natural frequencies $\omega_k = \Gamma \eta_k$ will become too disparate for the coupling to lock the oscillators to a common $\Omega$. So as $\Gamma$ increases, we eventually lack any stable phase-locked solution. 

This desynchronization transition occurs at what we call the {\it locking threshold}, at a parameter value given by the {\it critical value} of $\Gamma$ (alternatively, the critical width). Its calculation has been a focus of many prior studies of the Kuramoto model. Among these, a major point of variation has come from the choice of coupling topologies. The manner in which the oscillators are connected can have drastic effects on the behavior of the critical $\Gamma$, as has been demonstrated in work on complete graphs, one-dimensional chains and rings, two-dimensional square grids, three-dimensional cubic lattices, $d$-dimensional hypercubic lattices, random graphs, small-world and scale-free networks, and so on. For recent reviews, see Refs.~\cite{dorfler14, pikovsky15, rodrigues16}.

In this paper, we analyze a tractable situation where the dependence of the critical $\Gamma$ on topology, as opposed to dimension, can be well characterized. Namely, if we have a one-dimensional lattice of oscillators with nearest-neighbor coupling, how does the critical $\Gamma$ depend on the choice of boundary condition? If oscillators 1 and $N$ are coupled, we call this the {\it ring} topology and denote its locking threshold by $\Gamma_R$. Alternatively, if oscillators 1 and $N$ are not connected, we call this the {\it chain} topology and write its corresponding locking threshold  as $\Gamma_C$. 

Intuitively, one might expect a ring and a chain to have similar locking thresholds, especially when $N$ becomes large. After all, the two topologies differ only by a single edge. On the other hand, that single edge is responsible for a topological (and hence \emph{qualitative}, not merely quantitative) change in the lattice's connectivity structure. For that reason it could conceivably have a very potent effect. 

Although the setting of one-dimensional lattices may seem overly simplistic, it has the advantage that both rings and chains of oscillators have been studied extensively, using various techniques to analyze their dynamics and bifurcations \cite{cohen82, ermentrout84, ermentrout85, kopell86, sakaguchi87, sakaguchi88, strogatz88a, strogatz88b, wiley06, muruganandam08, elnashar09, ochab09, kogan09, lee09, giver11, tilles11, tilles13a,  tilles13b}.

The main question is this: If we have a chain and a ring subject to the same initial condition $\theta_k(0)=0$ and the same vector of base frequencies $\eta = (\eta_1, ..., \eta_N)$, what limits can be placed on the ratio $\Gamma_R/\Gamma_C$? In particular, must a ring always be ``more stable'' than a chain, leading to $\Gamma_R/\Gamma_C \geq 1$? 

\subsection{Telescopic Coupling}

In addition to variation in the connectivity structure, another source of variation in Kuramoto-like systems comes from altering the coupling function. For instance, we could replace the pure sine function in Eq.~\eqref{GeneralKura} with a more general periodic function.  As we will see, the following analysis allows for such a generalization, though at the cost of introducing a different type of special structure. 

To motivate this structure, let us look at the governing equation for an internal oscillator $k$ (meaning an oscillator with $1< k < N$) in a one-dimensional Kuramoto chain or ring:
\begin{alignat}{1} \label{StandardSine}
\dot\theta_k= \omega_k + \sin(\theta_{k-1} - \theta_k) + \sin(\theta_{k+1} - \theta_k).
\end{alignat}
Because sine is odd, Eq.~\eqref{StandardSine} can be rewritten as 
\begin{alignat}{1} \label{TelescopeSine}
\dot\theta_k= \omega_k + \sin(\theta_{k-1} - \theta_k) - \sin(\theta_{k} - \theta_{k+1}).
\end{alignat}
So if we want to generalize from sine to a more general function $f$, mathematically speaking we have two plausible choices: Either
\begin{equation}
\dot\theta_k =  \omega_k + f(\theta_{k-1} - \theta_k) + f(\theta_{k+1} - \theta_k), \label{StandardEq}
\end{equation}
or 
\begin{equation}
\dot\theta_k = \omega_k + f(\theta_{k-1} - \theta_k) - f(\theta_k - \theta_{k+1}). \label{TelescopicEq}
\end{equation}

Equation~\eqref{StandardEq} is a generalization of the Kuramoto model that has often been studied in the past, motivated by its physical and biological applications~\cite{ermentrout84, ermentrout85, kopell86, ostborn04}. However, we believe it is instructive to consider the alternative Eq.~\eqref{TelescopicEq} as well, and will devote most of our attention to it below. Where the distinction becomes important, we will say Eq.~\eqref{StandardEq} represents {\it standard coupling}, and Eq.~\eqref{TelescopicEq} represents {\it telescopic coupling}, thanks to some convenient cancellation properties it enjoys. We will restrict attention to continuously differentiable coupling functions $f$ that are $2\pi$-periodic, and will also demand that $f$ is nonconstant and has at least one zero. 

Although telescopic coupling is unconventional, it coincides with standard coupling when $f$ is an odd function, as commonly assumed in the physics literature. In that sense, telescopic and standard coupling schemes are on equal footing as generalizations of Kuramoto's sinusoidal coupling. Actually, considering the vast literature that focuses on pure sine coupling, even that special case remains of interest. Our results for telescopic coupling will include the traditional sine case while extending it to a new and wider family of models.  

One possible objection is that telescopic coupling injects a directionality to a chain or ring. To see this, note that swapping the oscillators and natural frequencies ``left to right'' ($j \to N-j+1$ for $j=1,\ldots, N$) changes the governing equations for telescopic coupling, but not for standard coupling. But such a directionality may be reasonable in some contexts. For example, there are a number of physical and biological systems which have been modeled as directed chains of oscillators, such as central pattern generators for the swimming rhythm of lamprey~\cite{cohen82, kopell86, cohen92, ren00}. 

A related point in favor of telescopic coupling is that it eases the analysis of oscillator arrays whose coupling functions $f$ lack odd symmetry. Although some results have been obtained for non-odd coupling on a chain ~\cite{kopell86, ostborn04, strogatz88a}, these are rare. Much of the existing research in this field has relied on the oddness of the coupling function and struggled otherwise. As we will see, telescopic coupling handles non-odd functions without difficulty. 

\section{Critical Width $\Gamma_C$ for a Chain}
To begin the analysis, we calculate the critical width $\Gamma_C$ above which the chain has no phase-locked solutions \cite{strogatz88a, strogatz88b}. After including the chain boundary terms, and assuming telescopic coupling as in Eq.~\eqref{TelescopicEq}, the dynamics are given by  
\begin{alignat*}{1}
\dot \theta_1 =& \omega_1 - f(\theta_1 - \theta_2), \\
\dot \theta_k =& \omega_k +f(\theta_{k-1}- \theta_k) - f(\theta_{k}- \theta_{k+1}), \text{for} \ 1 < k < N, \\
\dot \theta_N =& \omega_N +f(\theta_{N-1}- \theta_N).
\end{alignat*} 
By definition, for $\Gamma \leq \Gamma_C$, the system evolves to a stable locked state, and conversely, locking is impossible for $\Gamma > \Gamma_C$. So if we find a condition on the existence of a locked state, we get a condition on $\Gamma_C$. 

Recall that locking occurs when $\dot\theta_k \equiv \Omega$ for all $k$, for some $\Omega$. If we simply sum all $N$ of the differential equations above and then divide by $N$, we find  
\begin{alignat*}{1}
\Omega =  \frac{1}{N}\sum_{k=1}^N \omega_k,
\end{alignat*}
where we took advantage of telescoping nature of telescopic coupling. Let 
\begin{equation}
\bar\omega =  \frac{1}{N}\sum_{k=1}^N \omega_k. 
\end{equation}
This allows us to rewrite our condition for locking as 
\begin{alignat*}{1}
\omega_1 - \bar \omega =& f(\theta_1 - \theta_2), \\
\omega_k - \bar \omega =& -f(\theta_{k-1}- \theta_k) + f(\theta_{k}- \theta_{k+1}), 1<k<N, \\
\omega_N - \bar \omega =& -f(\theta_{N-1}- \theta_N).
\end{alignat*} 
Sum the first $k$ equations and telescope them to obtain 
\begin{alignat*}{1}
\sum_{j=1}^k \left(\omega_j - \bar \omega\right) = f(\theta_k - \theta_{k+1}).
\end{alignat*}
Let us define 
\begin{equation}
\phi_{k} = \theta_k - \theta_{k+1}
\end{equation} 
and 
\begin{equation}
D_k = \sum_{j=1}^k (\eta_j - \bar \eta)
\end{equation} 
for $k = 1, \ldots, N-1$. This yields 
\begin{alignat}{1} \label{ChainEq}
f(\phi_k) = \Gamma D_k,
\end{alignat}
which is an exact condition on finding a locked state in the chain topology. In particular this means that $\Gamma_C$ corresponds to the supremum of all $\Gamma$'s where the above equation is satisfied and the solution $\phi = (\phi_1, .., \phi_{N-1})$ is stable. This condition is equivalent to one  found previously for sine coupling~\cite{strogatz88a, strogatz88b}. 

\subsection{Existence and Stability of the Locked State}
Next we check that condition \eqref{ChainEq} is satisfiable for the class of $f$ under consideration. Our biggest demand on $f$ was that it be continuously differentiable and periodic. Continuous periodic real-valued functions are bounded and attain their maximum and minimum,  so we know that both $f$ and $f'$ attain their upper and lower bounds. Let us define the bounds $f_u := \max_x f(x),$ and $f_l := \min_x f(x).$ We also requested that $f$ be non-constant and cross zero, so $f_u > 0 > f_l.$

Given a particular realization of $\eta_k$'s, we can define $D_u := \max(0, \max_k(D_k))$ and $D_l := \min(0, \min_k(D_k))$. So $D_u$ represents the largest positive value of $D_k$ if it exists and 0 otherwise, with $D_l$ similarly defined for negative values, enforcing $D_l \leq 0 \leq D_u$. Therefore, we know that all locked states disappear at a critical point of 
\begin{equation} \label{ChainCrit} \Gamma_C = \min(f_u/D_u, f_l/D_l ).\end{equation}
We formally take $1/0 = \infty$; notice that $\Gamma_C = \infty$ if and only if $D_k =0$ for all $k$, which is only possible if all the $\eta_k$ are identical. Also note that since $f_u$ and $f_l$ represent global bounds on $f$, then no equilibrium {\it at all} can exist when $\Gamma > \Gamma_C.$ However, for $\Gamma <\Gamma_C$ we can always find a set of $\phi_k$ that will satisfy the prior equations. This makes $\Gamma_C$ the true point between a locked state existing and disappearing. 

\begin{figure}
\includegraphics[width = 0.5\textwidth]{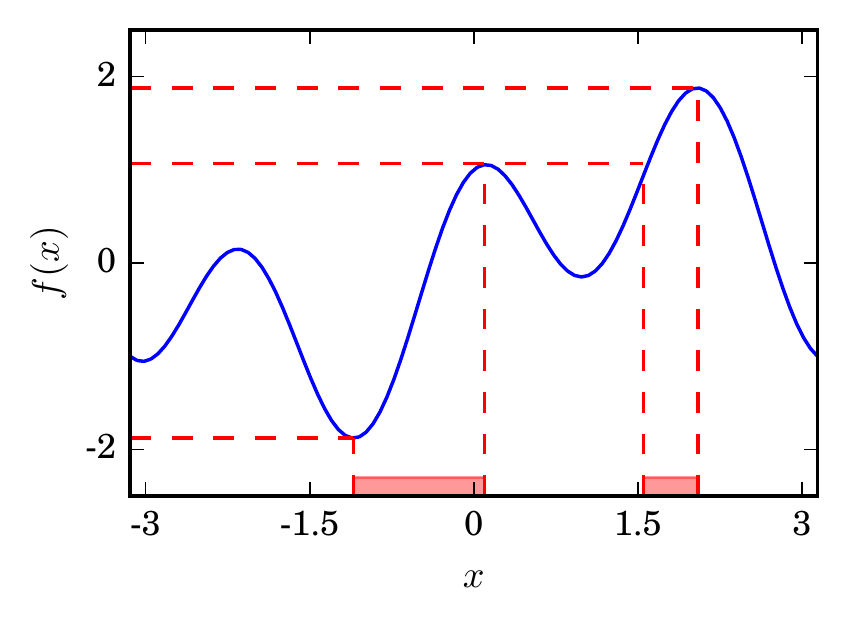}
\caption{Example showing a choice of $\Lambda$ for a specific coupling function $f(x) = \sin(x)+\cos(3x)$, as indicated by the shaded region on the $x$-axis. The dashed lines illustrate how $\mbox{image}(f|_\Lambda)$ $= \mbox{image}(f)$, up to global extrema, while having $f'|_\Lambda(x) > 0 $. }
\label{LambdaPlot}
\end{figure}

However, knowing a phase-locked state exists does not ensure it is stable. Fortunately, it is not hard to show if $\Gamma < \Gamma_C$, then a stable locked state exists. For any $y \in (f_l, f_u)$ there exists some point $x$ where $f(x) = y$ and $f'(x) >0$; otherwise $f$ could never climb from $y$ to $f_u$. Moreover, since $f'$ is bounded, there are only finitely many $x$ which could work for each $y$ in the bounded domain $(-\pi, \pi]$. And since $f'$ is continuous, then $f'$ will be positive in a neighborhood of $x$, so our point selection can take advantage of this. Ergo there exists some open set $\Lambda$ where $f$ restricted to $\Lambda$ always has positive derivative and is surjective onto $(f_l, f_u).$  For a visual example, see Fig.~\ref{LambdaPlot}.

Returning to our original question, if $\Gamma < \Gamma_C,$ we can select a set of $\phi_k$ out of $\Lambda$ in a well-defined way, where $f'(\phi_k) >0$ and $f(\phi_k) = \Gamma D_k$. A theorem of Ermentrout ~\cite{ermentrout92} then guarantees that such a solution is asymptotically stable. Therefore, Eq.~\eqref{ChainCrit} really does define $\Gamma_C$, below which at least one stable locked state exists and above which none do.

\section{An upper bound on $\Gamma_R/\Gamma_C$}
The next step is to obtain an upper bound on $\Gamma_R$, the locking threshold for a ring. Although the interior of a chain looks the same as a ring, they differ at the boundary terms, as seen in the following equations:  
\begin{alignat*}{1}
\dot \theta_1 =& \omega_1 + f(\theta_N - \theta_1) - f(\theta_1 - \theta_2), \\
\dot \theta_k =& \omega_k +f(\theta_{k-1}- \theta_k) - f(\theta_{k}- \theta_{k+1}), 1<k<N, \\
\dot \theta_N =& \omega_N +f(\theta_{N-1}- \theta_N) - f(\theta_N - \theta_1).
\end{alignat*} 
Nevertheless, several steps in the following argument will be the same as for the chain. For example, locked states still satisfy  $\dot\theta_k = \Omega$ for some $\Omega$, and we can still telescope the equations, yielding $\Omega = \bar\omega$ again.  Similarly,
\begin{alignat*}{1}
\omega_1 - \bar \omega =& -f(\theta_N - \theta_1) + f(\theta_1 - \theta_2), \\
\omega_k - \bar \omega =& -f(\theta_{k-1}- \theta_k) + f(\theta_{k-}- \theta_{k+1}), 1<k<N, \\
\omega_N - \bar \omega =& -f(\theta_{N-1}- \theta_N) + f(\theta_N - \theta_1)
\end{alignat*} 
which can be telescoped into
\begin{equation}
\Gamma D_k = f(\phi_k) - f(\theta_N - \theta_1). \nonumber
\end{equation}
Here, $D_k$ and $\phi_k$ are defined exactly as in the last section. Hence, if we put the same choice of $\eta$'s on a ring and a chain, they would have the same vectors $D = (D_1, ..., D_{N-1})$. Also notice that $-\sum_{j=1}^{N-1} \phi_j = \sum_{j=1}^{N-1} (\theta_j - \theta_{j-1}) = \theta_N - \theta_1.$ Therefore, we can write
\begin{equation}\label{RingEq} 
\Gamma D_k = f(\phi_k) - f\left(-\sum_{j=1}^{N-1} \phi_j\right).
\end{equation}

Equations like this have been found before for the special case of sine coupling~\cite{ochab09, tilles11, tilles13b}. Although Eq.~\eqref{RingEq} has a compact form, demonstrating that solutions to it exist and calculating them explicitly is a difficult endeavor; hence our more modest goal is to establish a bound on $\Gamma_R$. 

Let us define $f_u, f_l, D_u$, and $D_l$ as before. Since $f_l < 0 < f_u$ represent global extrema, the ring can have a locked state only if $f_l - f_u \leq \Gamma D_k \leq f_u - f_l$ for all $k$. This yields  
\begin{equation}\label{RingCrit}
\Gamma_R \leq \min\left(\frac{f_u - f_l}{D_u}, \frac{f_l - f_u}{D_l} \right).
\end{equation}
Because we have been careful to use the same $D_k$ here as in the chain case, Eq.~\eqref{RingCrit} can be directly compared to Eq.~\eqref{ChainCrit} to give the bound  
\begin{alignat*}{1}
\Gamma_R/\Gamma_C \leq& \min\left(\frac{f_u - f_l}{D_u}, \frac{f_l - f_u}{D_l} \right) / \min\left(\frac{f_u}{D_u},\frac{f_l}{D_l} \right) \\
=& \min\left(\frac{f_u - f_l}{D_u}, \frac{f_l - f_u}{D_l} \right)\max\left(\frac{D_u}{f_u}, \frac{D_l}{f_l} \right). \notag 
\end{alignat*}
If $D_u/f_u > D_l/f_l$, then we have that 
\begin{alignat*}{1}
\Gamma_R/\Gamma_C \leq& \min\left(\frac{f_u - f_l}{D_u}, \frac{f_l - f_u}{D_l} \right) \left( \frac{D_u}{f_u} \right) \\
\leq & \left(\frac{f_u - f_l}{D_u}\right) \left( \frac{D_u}{f_u}\right) \\
=& 1 + \left|\frac{f_l}{f_u} \right|.
\end{alignat*}
If otherwise $D_u/f_u < D_u/f_l$, then
\begin{alignat*}{1}
\Gamma_R/\Gamma_C \leq& \min\left(\frac{f_u - f_l}{D_u}, \frac{f_l - f_u}{D_l} \right) \left( \frac{D_l}{f_l} \right) \\
\leq& \left(\frac{f_l - f_u}{D_l}\right) \left( \frac{D_l}{f_l}\right)\\
=& 1 + \left|\frac{f_u}{f_l} \right|.
\end{alignat*}
Together these facts imply  
\begin{equation} \label{RatioBound}
\Gamma_R/\Gamma_C \leq 1 + \max\left(\left| \frac{f_l}{f_u} \right|, \left|\frac{f_u}{f_l} \right| \right).
\end{equation}

Thus we have found a rigorous upper bound on the ``advantage'' of a ring over a chain, with the bound depending exclusively on the shape of the coupling function $f$. Also notice that the arguments of the $\max$ function are a nonnegative real number and its reciprocal, so this upper bound is always at least 2. 

In fact, if $f$ is odd and  $N=2$ then $\Gamma_R = 2 \Gamma_C,$ so Eq.~\eqref{RatioBound} is tight in certain cases. We will discuss the linear stability of states on the ring later, but since we need a solution to exist before it can be stable, the bound is valid.

\section{Upper and lower bounds}

\begin{figure*}
\includegraphics[width = \textwidth]{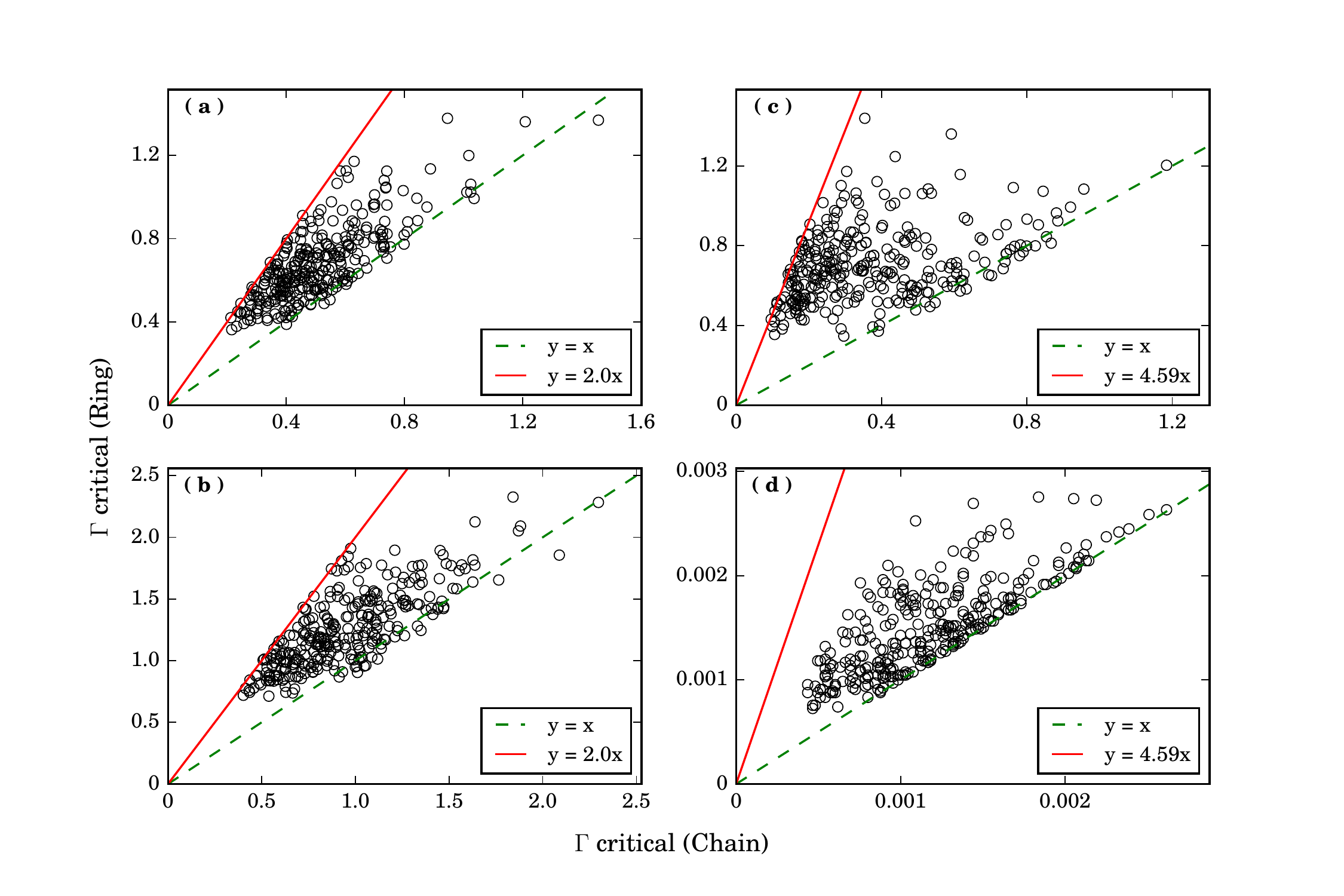}
\caption{Scatterplot comparing the critical width $\Gamma$ for a chain and a ring of $N=25$ oscillators for a variety of coupling schemes and random realizations of the natural frequencies. For each data point, the corresponding ring and chain were matched, meaning that both were  subject to the same initial conditions and natural frequencies. Initial phases $\theta_k(0)$ were chosen to be identically zero, and natural frequencies $\eta_k$ were drawn at random from a uniform distribution on $[-1,+1]$. Lines here represent a 1:1 ratio (dashed green line) and our theoretically predicted upper bound (solid red line defined by Eq.~\eqref{RatioBound}). Panel (a) has $f(x) = \sin(x)$ and panel (b) has $f(x) = \sin(x)+\cos(3x)$, both under the telescopic coupling scheme of Eq.~\eqref{TelescopicEq}. Both the panels (c) and (d) on the right have $f(x) = \sin(x +0.6) - \sin(0.6)$. However (c) uses telescopic coupling, whereas (d) follows the standard coupling equations of \eqref{StandardEq}. Notice that the upper bound is always obeyed, but some data points lie below the lower dashed line, showing that it is not a strict bound. Values of $\Gamma$ were estimated via a bisection technique combined with numerical integration, using a fourth-order Runge-Kutta method with a timestep of 0.125, a transient time between $5\ \times 10^2$ and $2 \times10^3$ time units, and observation times of $5 \times 10^2$. }
\label{WedgePlot}
\end{figure*}

Now that we have the upper bound \eqref{RatioBound} on the ratio of the critical widths, it is natural to want to check how sharp it is. The results shown in Fig.~\ref{WedgePlot} do exactly that. We generate many different realizations of the base frequency vectors $\eta$, and then plot the numerically obtained $\Gamma_C$ and $\Gamma_R$ on a scatterplot, and draw a solid line to denote our predicted boundary \eqref{RatioBound}. We first test an odd coupling function, namely sine;  then we test  several non-odd coupling functions. For the regimes being tested, a lot of points congregate at our upper bound, but as expected, none actually trespass it.

But what about the lower dashed line representing $\Gamma_C = \Gamma_R$? It is tempting to think that this line should also be respected; after all, a ring has an additional coupling connection, and it has no free ends. With this extra edge to provide more coupling between the oscillators, one intuitively expects that a ring should always lock more easily than a chain. Moreover, the difference in boundary conditions means that the ring permits topologically twisted  states ~\cite{ermentrout85, wiley06} that would be impossible for the chain. This too would naively suggest that the ring is always more susceptible to locking than the chain is. 

However, Fig.~\ref{WedgePlot} indicates that some cases lie below the dashed line. In such cases the chain locks when its matched ring does not. Apparently the naive intuition above is wrong. We now confirm this surprising result by constructing a counterexample.

\subsection{Counterexample to $\Gamma_C \leq \Gamma_R$}

In fact, the critical width of a chain is \emph{not} always less than that of a matched ring. Here is a counterexample. Say we have $N = 4$, $f = \sin$, and we have obtained a realization of $\eta$'s such that $D = (+1, -1, -1)^T$. Then Eq.~\eqref{ChainCrit} immediately implies that $\Gamma_C = 1$, and we can satisfy this system with $\phi_1 = -\phi_2 = -\phi_3 = \pi/2$. 

Now consider what the corresponding locked state would be for the ring.  By assumption, such a state must exist; if $\Gamma_C \leq \Gamma_R$ is true, we should be able to produce a locked solution to the ring equations~\eqref{RingEq} with $\Gamma=1$. Such a solution would then satisfy the following system:
\begin{alignat}{1} \label{Countereqs}
\sin(\phi_1) + \sin(\phi_1+\phi_2+\phi_3) &= +1, \\
\sin(\phi_2) + \sin(\phi_1+\phi_2+\phi_3) &= -1, \\
\sin(\phi_3) + \sin(\phi_1+\phi_2+\phi_3) &= -1.
\end{alignat}
Notice that if we subtract the second or third equation from the first, we get
\begin{alignat*}{1} 
\sin(\phi_1) - \sin(\phi_2) &= +2, \\
\sin(\phi_1) - \sin(\phi_3) &= +2.
\end{alignat*}

From here, we realize we have no choice. It must be that $\phi_1 = \pi/2$ and $\phi_2 = \phi_3 = -\pi/2,$ which yields the desired contradiction, since it gives 
\begin{equation} 
\sin(\phi_1) + \sin(\phi_1+\phi_2+\phi_3) = \sin(\pi/2) + \sin(-\pi/2) = 0, \nonumber 
\end{equation}
which violates Eq.~(14). 

The contradiction shows that even though we have a locked state for a chain, none exists for the ring. So sometimes $\Gamma_C \not\leq \Gamma_R$. Because this counterexample uses the sine function, it works for both the standard and telescopic coupling models. 

Unfortunately, trying to come up with a genuine lower bound for $\Gamma_R/\Gamma_C$ is surprisingly involved, given that such shenanigans can be found in the small-$N$ cases.

\section{Asymptotic existence}

Although small $N$ is problematic, the large-$N$ regime is more tractable. Let us fix $N$ to be large but finite, and choose some realization of $\eta$. If we start with $\Gamma < \Gamma_C$, we are guaranteed a phase-locked solution $\phi^{(C)}$ to the chain equation \eqref{ChainEq}, satisfying $f\left(\phi_{k}^{(C)}\right) = \Gamma D_k$ for all $k = 1, \dots, N-1$. Moreover, we are guaranteed to be able to choose these $\phi_k$ from the set $\Lambda$ as defined earlier.

We seek to construct an approximate phase-locked solution to the ring based on this chain solution. The coupling function $f$ is $2\pi$-periodic, so let us define 
\[\Psi := \left( \sum_{j=1}^{N-1} \phi_{k}^{(C)}\right)\mbox{mod}{2\pi}. \] 
Thus $0 \leq \Psi < 2\pi$.  Since we insisted that $f$ cross zero and be both nonconstant and periodic, there exists some point $x_0 \in (-\pi,\pi]$ such that $f(x_0) = 0$ and $f'(x_0) > 0$. We can then define $ \phi_{k}^{(R)} := \phi_{k}^{(C)} - (x_0 + \Psi)/(N-1)$ as a value close to $\phi_{k}^{(C)}$. This will represent our attempted solution to the ring equation \eqref{RingEq}. 

First notice that 
\begin{alignat*}{1}
f\left(-\sum_{j=1}^{N-1} \phi_{j}^{(R)}\right) &= f\left(-\sum_{j=1}^{N-1} \left(\phi_{k}^{(C)} - \frac{x_0 + \Psi}{N-1} \right) \right) \\
&= f\left( -\sum_{j=1}^{N-1} \left(\phi_{k}^{(C)}\right) + x_0 + \Psi \right)\\
&= f(x_0)\\
&= 0.
\end{alignat*}
And so we find 

\begin{alignat*}{1}
&f\left(\phi_{k}^{(R)}\right) - f\left(-\sum_{j=1}^{N-1} \phi_{j}^{(R)}\right) \\
&= f\left(\phi_{k}^{(C)} - \frac{x_0 + \Psi}{N-1}\right).
\end{alignat*}
But recall that $f$ is continuously differentiable, so there is some finite upper bound on the derivative $f_u' = \max_x |f'(x)|$. In other words, for any $x$ and $\delta$, then $|f(x) - f(x+\delta)| < f_u' \delta.$ Therefore, 
\begin{equation*}
\left| f\left(\phi_{k}^{(C)} - \frac{x_0 + \Psi}{N-1} \right) - f\left(\phi_{k}^{(C)}\right) \right| < f_u' \frac{|\Psi+x_0|}{N-1},
\end{equation*}
which implies
\begin{alignat*}{1}
f\left(\phi_{k}^{(R)}\right) - f\left(-\sum_{j=1}^{N-1} \phi_{j}^{(R)}\right) =\Gamma  D_k + O(N^{-1}).
\end{alignat*}
So $\phi^{(R)}$ is an approximate solution to the ring equations that becomes exact as $N$ approaches infinity. 

\begin{figure}
\centering
\includegraphics[width = 0.5\textwidth]{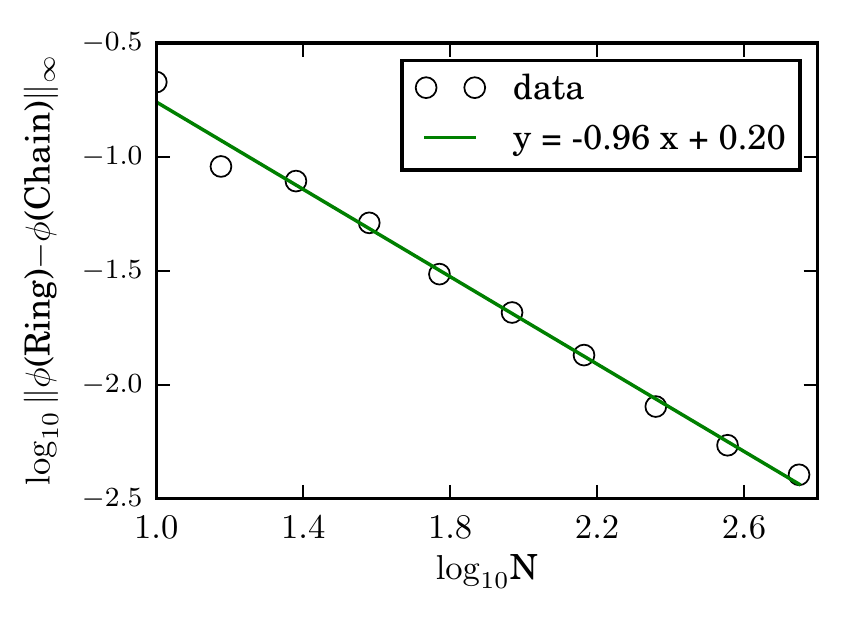}
\caption{A plot showing the log of the separation (as measured by the infinity norm) between the vector of $\phi_k$'s for a chain and the same vector for a ring, given that they are subject to the same natural frequencies $\eta_k$, which were randomly drawn from a uniform distribution on $[-1,+1]$. We first calculate the locked solution for a chain, using an initial condition of all zeros. Then we use the final result of that calculation as the initial condition of the ring to allow for direct comparison. A coupling function $f(x) = -\sin(x)$ was used. The straight line shows the best linear fit to the log-log plot, indicating that we are seeing a decay comparable to $O(N^{-1})$. Values of the phases were computed by numerical integration,  using a fourth-order Runge-Kutta method with a timestep of 0.125, a transient time between $5 \times10^2$ and $10^6$ time units and observation times of $10^3$.}
\label{NormPlot}
\end{figure}

Figure~\ref{NormPlot} shows the convergence of this approximate solution for the ring to that for the chain. We numerically construct pairs of solutions that get closer as $N$ gets large. This all makes sense, since an infinitely long chain should be identical to an infinitely long ring. 

Concerning stability, remember that the set $\Lambda$ is open, so for any $x\in \Lambda$, then for sufficiently small $\delta$ then $x+\delta \in \Lambda$. So this ring solution $\phi^{(R)}$ also lies entirely in $\Lambda$ for large enough $N$. This is almost enough to cite Ermentrout and establish the stability of this solution~\cite{ermentrout92}. However, we have an additional phase difference in our dynamics, $\theta_N - \theta_1 = -\sum_{j=1}^{N-1} \phi_{k}.$ In our proposed solution this quantity is sent to $\Psi+x_0$, which by construction has $f'(\Psi + x_0) = f'(x_0) > 0$, and so stability is secured. 

To summarize, if we have a stable locked solution to the chain of oscillators for large $N$, then there is a nearby stable locked solution for the ring of oscillators. Hence, the naive lower bound $\Gamma_R \geq \Gamma_C$ is valid in the asymptotic case $N \gg 1$.

\subsection{Partial Results for Standard Coupling}
Our prior argument relied very little on telescopic coupling. In fact, we can use a similar method to show an equivalent result using the standard coupling~\eqref{StandardEq} instead of telescopic coupling~\eqref{TelescopicEq}. However, this requires the additional constraint of $x_0 = 0$ or $\pi$ (where $f(x_0) = 0$ and $f'(x_0) > 0$). 

To derive the relevant results, suppose that we have some set of $\phi_{k}^{(C)} = \theta_{k} - \theta_{k+1}$ which satisfy the standard coupling equations for the chain and are locked at $\dot\theta_k \equiv \Omega$. Then 
\begin{alignat*}{2}
& \Omega = \omega_1 + f\left(-\phi_{1}^{(C)}\right), \\
& \Omega = \omega_k +f\left(\phi_{k-1}^{(C)}\right) + f\left(-\phi^{(C)}_{k} \right), 1< k <N, \\
& \Omega = \omega_N + f\left(\phi^{(C)}_{N-1}\right).
\end{alignat*} 
Using the fact that $\theta_N - \theta_0 = -\sum_{j=1}^{N-1} \phi_{k}$, the condition for locking on a ring becomes 
\begin{alignat*}{2}
& \Omega = \omega_1 + f\left(-\phi_{1}^{(R)}\right) + f\left(-\sum_{j=1}^{N-1} \phi_{k}^{(R)} \right), \\
& \Omega = \omega_k +f\left(\phi_{k-1}^{(R)}\right) + f\left(-\phi^{(R)}_{k} \right), 1<k <N, \\
& \Omega = \omega_N + f\left(\phi^{(R)}_{N-1}\right) + f\left(\sum_{j=1}^{N-1} \phi_{k}^{(R)} \right).
\end{alignat*} 
If we try plugging in $ \phi_{k}^{(R)} := \phi_{k}^{(C)} - (\Psi+x_0)/(N-1)$, with $\Psi$ defined the same as before, then the sum terms will evaluate to $x_0$ modulo $2\pi$. If we use the  continuity arguments from before for $1 < k <N$, then
\begin{alignat*}{1}
& \omega_k +f\left(\phi_{k-1}^{(C)}- \frac{\Psi + x_0}{N-1} \right) + f\left(-\phi^{(C)}_{k} + \frac{\Psi + x_0}{N-1} \right) \\
&= \Omega + O(N^{-1}).
\end{alignat*}
For $k = 1$, then 
\begin{alignat*}{2}
& \omega_1 + f\left(-\phi_{1}^{(R)}\right) + f\left(-\sum_{j=1}^{N-1} \phi_{k}^{(R)} \right) \\
&= \omega_1 + f\left(-\phi_{1}^{(C)} +\frac{\Psi + x_0}{N-1} \right)+ f(x_0) \\
&= \Omega + O(N^{-1}),
\end{alignat*}
and for $k=N$, then
\begin{alignat*}{1}
& \omega_N + f\left(\phi^{(R)}_{N} \right) + f\left(\sum_{j=1}^{N-1} \phi_{k}^{(R)} \right)\\
&= \omega_N + f\left(\phi_{N}^{(C)} - \frac{\Psi + x_0}{N-1} \right)+ f(-x_0) \\
&= \Omega + O(N^{-1}).
\end{alignat*} 
Hence, as $N \rightarrow \infty$  this solution becomes exact. The reason we restricted $x_0$ was because we wanted $f(x_0) = 0 = f(-x_0)$, which was only guaranteed if $x_0 = -x_0$ mod $2\pi$. 

By periodicity and continuity, if $f_l < f(\phi_k) < f_u$, then there is always some $\phi_k'$ such that $f(\phi_k') = f(\phi_k)$ and $f'(\phi_k') > 0$. So without loss of generality, if we had a chain solution $\phi_{k}^{(R)}$, we could pick another solution where all the phase differences have positive slope in $f$. And for sufficiently large $N$, the same would hold true for the ring, since we are perturbing only slightly and we already assumed $f'(x_0)>0$ for the boundary term. 

Therefore, given any existing locked solution for the chain with standard coupling (even an unstable solution), this argument guarantees the existence of a \emph{stable} locked solution for the chain and a stable approximate solution for the ring, also with standard coupling. This means the large-$N$ limit gives $\Gamma_R \geq \Gamma_C$ for standard coupling, just as it did for telescopic coupling. But unlike the more convenient case of telescopic coupling, we can no longer construct a locked solution in the first place nor can we put clean upper bounds on $\Gamma_R$ or $\Gamma_C$.

\section{Summary and future directions}

\begin{figure}
\centering
\includegraphics[width = 0.5\textwidth]{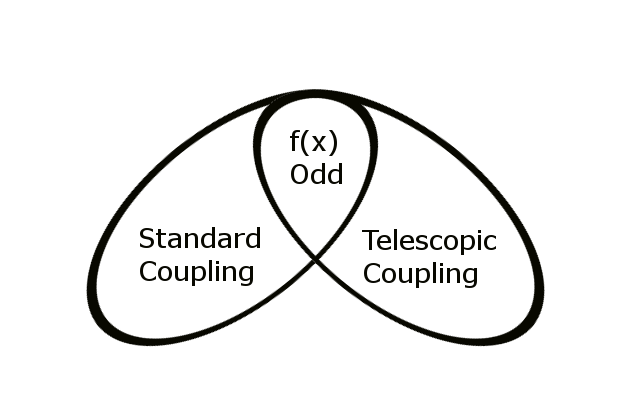}
\caption{Schematic illustration of the relationship between standard coupling, telescopic coupling, and their agreement for odd coupling functions $f$.}
\label{ModelSchematic}
\end{figure}

Our main results have put a limit on the relative behavioral difference between a ring and a chain of phase oscillators. As noted in the introduction, it is generically hard to predict the conditions on synchronization. As we have shown, even a single additional connection can cause a  \emph{doubling} of the locking threshold, emphasizing its sensitivity to topology. However, simply putting limits on the synchronization criterion is often good enough for practical purposes. This is especially true for our particular comparison, since an analytic criterion is exactly known for a chain, but no equivalent has been demonstrated for a ring. 

Our analysis was facilitated by the introduction of the telescopic coupling scheme \eqref{TelescopicEq}. Thanks to its convenient analytic properties, a large collection of different results which typically require sine, odd, or some other heavily restricted coupling function have been generalized to a new family of $f$'s. And as we noted earlier, telescopic coupling \eqref{TelescopicEq} and standard coupling \eqref{StandardEq} have equally legitimate mathematical claims to being a generalization of the sine-based Kuramoto model. Moreover, as illustrated by Fig.~\ref{ModelSchematic}, these two coupling schemes exactly overlap in the case of odd $f$. 

Regarding future directions, one possibility is to ask whether the results generalize to higher dimensions. Telescopic coupling introduces a directionality to a one-dimensional chain. The natural extension to higher-dimensional lattices would be to introduce a directionality along each axis. It is not hard to do such a thing, and when we do, we extend the same cancellation properties enjoyed by odd $f$ to generic $f$. What results might come of this? 

Although we distinguished between the two coupling schemes  \eqref{StandardEq} and \eqref{TelescopicEq}, we have not made much effort to connect their behaviors for $f$ non-odd. But hints of a connection are present. For example, the two plots on the right side of Fig.~\ref{WedgePlot} both seemingly obey our predictions, even though only one of them used our preferred telescopic coupling. The other came from standard coupling, about which we were unable to make comparably strong statements. Unfortunately, we made little progress in computing exact relationships between the two schemes. Given that the telescoping scheme is much easier to work with, any general connection could potentially shed a lot of light on the standard case.

Finally, the fundamental question of this paper could be generalized in an ambitious way. Given a network of oscillators on a connectivity graph $G$, how does the locking threshold $\Gamma$ change with the addition or removal of a single edge? We chose both the graph and the edge very carefully in this paper, but some of our logic might still be relevant to this larger problem. Considering the close relationship between the power grid and the Kuramoto model, this question might bear on current issues of power grid resilience~\cite{dorfler13, dorfler12, motter13, wang09, kinney05, albert04}.

This research was supported by a Sloan Fellowship and NSF Graduate Research Fellowship grant DGE-1650441 to Bertrand Ottino-L\"{o}ffler in the Center for Applied Mathematics at Cornell, as well as by NSF grants DMS-1513179 and CCF-1522054 to Steven Strogatz. 


\end{document}